\documentclass[12pt]{article}
\usepackage{amsmath,amssymb,amsthm}

\newtheorem{theorem}{Theorem}
\newtheorem{question}{Question}
\newtheorem{corollary}{Corollary}

\begin{document}
\begin{center}
\Large
Digit Reversal Without Apology
\end{center}

\begin{flushright}
Lara Pudwell  \\
Rutgers University \\
Piscataway, NJ 08854\\
\texttt{lpudwell@math.rutgers.edu}

\vspace{2 mm}

\end{flushright}

In \emph{A Mathematician's Apology}~\cite{gH93} G. H. Hardy states, ``$8712$ and $9801$ are the only four-figure numbers which are integral multiples of their reversals"; and, he further comments that ``this is not a serious theorem, as it is not capable of any significant generalization."

However, Hardy's comment may have been short-sighted.  In 1966, A. Sutcliffe~\cite{aS66} expanded this obscure fact about reversals.  Instead of restricting his study to base $10$ integers and their reversals, Sutcliffe generalized the problem to study all integer solutions of 
$$k(a_hn^{h}+a_{h-1}n^{h-1}+ \cdots + a_1 n+a_0) = a_0 n^{h} + a_1 n^{h-1}+ \cdots + a_{h-1}n + a_h$$ with $n \geq 2$, $1<k<n$, $0 \leq a_i \leq n-1 \mbox{ for all } i$, $a_0 \neq 0$, $a_h \neq 0$.  We shall refer to such an integer $a_0...a_h$ as an \emph{$(h+1)$-digit solution for $n$} and write $k(a_h,a_{h-1},...,a_1,a_0)_n = (a_0,a_1,...,a_{h-1},a_h)_n$.  For example, $8712$ and $9801$ are $4$-digit solutions in base $n=10$ for $k=4$ and $k=9$ respectively.  After characterizing all $2$-digit solutions for fixed $n$ and generating parametric solutions for higher digit solutions, Sutcliffe left the following open question: Is there any base $n$ for which there is a $3$-digit solution but no $2$-digit solution?

Two years later T. J. Kaczynski \footnote{Better known for other work.}~\cite{tK68} answered Sutcliffe's question in the negative.  His elegant proof showed that if there exists a $3$-digit solution for $n$, then deleting the middle digit gives a $2$-digit solution for $n$.  Together with Sutcliffe's work, this proved that there exists a $2$-digit solution for $n$ if and only if there exists a $3$-digit solution for $n$.

Given the nice correspondence between $2$- and $3$-digit solutions described by Sutcliffe and Kaczynski, it is natural to ask if there exists such a correspondence for higher digit solutions.  In this paper, we will explore the relationship between $4$- and $5$-digit solutions.  Unfortunately, there is not a bijection between these solutions, but there is a nice family of $4$- and $5$- digit solutions which have a natural one-to-one correspondence.

A second extension of Sutcliffe and Kaczynski's results is to ask, ``Is there any value of $n$ for which there is a $5$-digit solution but no $4$-digit solution?"  We will answer this question in the negative; and, furthermore, we will show that there exist $4$- and $5$-digit solutions for every $n \geq 3$.

\subsection*{An attempt at generalization}\label{S:generalize}

In the case of $3$-digit solutions, Kaczynski proved that if $n+1$ is prime and $k(a,b,c)_n=(c,b,a)_n$ is a $3$-digit solution for $n$, then $k(a,b)_n = (b,a)_n$ is a $2$-digit solution.  Thus, we consider the following:

\begin{question}
Let $k(a,b,c,d,e)_n = (e,d,c,b,a)_n$ be a 5-digit solution for $n$.  If $n+1$ is prime, then is $k(a,b,d,e)_n = (e,d,b,a)_n$ a 4-digit solution for $n$?
\end{question}

First, following  Kaczynski, let $p=n+1$.  We have 
\begin{equation}\label{E:5sol}
k(an^4+bn^3+cn^2+dn+e) = en^4+dn^3+cn^2+bn+a.
\end{equation}
Reducing this equation modulo $p$, we obtain
$$k(a-b+c-d+e) \equiv e-d+c-b+a = a-b+c-d+e  \mbox{ mod } p.$$
Thus, $(k-1)(a-b+c-d+e) \equiv 0 \mbox{ mod } p$, and 
\begin{equation}\label{E:divisor}
p \mid (k-1)(a-b+c-d+e).
\end{equation}
If $p \mid (k-1)$, then $k-1 \geq p$, which is impossible because $k<n$.  Therefore, $p \mid (a-b+c-d+e)$.  But $-2p < -2n < a-b+c-d+e < 3n <3p$, so there are four possibilities:
\\(i) $a-b+c-d+e=-p$,
\\(ii) $a-b+c-d+e=0$,
\\(iii) $a-b+c-d+e=p$,
\\(iv) $a-b+c-d+e=2p$.
\\

Write $a-b+c-d+e=fp$, where $f \in \{-1,0,1,2\}$.  Substituting $c=-a+b+d-e+fp$ into equation \ref{E:5sol} gives:
$$k[n^2(n^2-1)a+n^2(n+1)b+fpn^2+n(n+1)d-(n^2-1)e]$$
$$= n^2(n^2-1)e+n^2(n+1)d+fpn^2+n(n+1)b-(n^2-1)a.$$

After substituting for $p$, dividing by $n+1$, and rearranging, one sees that $k[an^3+(b-a+f)n^2+(d-e)n+e]=en^3+(d-e+f)n^2+(b-a)n+a$.  Indeed, this is a $4$-digit solution for $n$ if $f=0$, $b-a\geq 0$, and $d-e\geq 0$, but not necessarily a $4$-digit solution of the form conjectured in Question 1.

As in Kaczynski's proof for $2$- and $3$-digit solutions, it would be ideal if three of the four possible values for $f$ lead to contradictions and the fourth leads to a ``nice" pairing of $4$- and $5$-digit solutions.  Unlike Kaczynski, we now have the added advantage of exploring these cases with computer programs such as Maple. Experimental evidence suggests that the cases $f=-1$ and $f=2$ are impossible.  The cases $f=0$ and $f=1$ are discussed below.

\subsection*{A counterexample}\label{S:counter}

Unfortunately, Kaczynski's proof does not completely generalize to higher digit solutions.  Most $5$-digit solutions do, in fact, yield $4$-digit solutions in the manner described in Question 1, but for sufficiently large $n$ there are examples where $(a,b,c,d,e)_n$ is a $5$-digit solution but $(a,b,d,e)_n$ is not a $4$-digit solution.

A computer search shows that the smallest such counterexamples appear when $n=22$:  
$$7(2,8,3,13,16)_{22} = (16,13,3,8,2)_{22}, 3(2,16,11,5,8)_{22} = (8,5,11,16,2)_{22}.$$  However, there is no integer $k$ for which $k(2,8,13,16)_{22} = (16,13,8,2)_{22}$ or $k(2,16,5,8)_{22} = (8,5,16,2)_{22}$. Note that $-2+8+13-16 = 3$ and $-2+16+5-8=11$; that is, both of these counterexamples to Question 1 occur when $f=0$.  The next smallest counterexamples are $$3(3,22,15,7,11)_{30}=(11,7,15,22,3)_{30}, 8(2,13,8,16,9)_{30}=(9,16,8,13,2),$$ which occur when $f=0$ and $n=30$.

\subsection*{A family of $4$- and $5$-digit solutions}\label{S:fam}

Although Kaczynski's proof does not generalize entirely, there exists a family of $5$-digit solutions when $f=1$ that has a nice structure.

\begin{theorem}Fix $n \geq 2$ and $a>0$.  Then
$$k(a,a-1,n-1,n-a-1,n-a)_n = (n-a,n-a-1,n-1,a-1,a)_n$$
 is a $5$-digit solution for $n$ if and only if $a \mid (n-a)$.
\end{theorem}

\begin{proof}  We have $$\frac{(n-a)n^4+(n-a-1)n^3+(n-1)n^2+(a-1)n+a}{an^4+(a-1)n^3+(n-1)n^2+(n-a-1)n+(n-a)}$$
$$= \frac{(n-a)(n^4+n^3-n-1)}{a(n^4+n^3-n-1)} = \frac{n-a}{a},$$ and the result is clear.
\end{proof}

Notice that $$(-a+(a-1))+((n-a-1)-(n-a))+p = -1+-1+(n+1) = n-1.$$
That is, this family of solutions occurs when $f=1$.  Moreover, this family follows the pattern described in Question 1; that is, for each $5$-digit solution described in Theorem 1, deleting its middle digit gives a $4$-digit solution.

\begin{theorem} If 
$$k(a,a-1,n-1,n-a-1,n-a)_n = (n-a,n-a-1,n-1,a-1,a)_n$$
 is a $5$-digit solution for $n$, then
$$k(a,a-1,n-a-1,n-a)_n = (n-a,n-a-1,a-1,a)_n$$
is a $4$-digit solution for $n$.
\end{theorem}

\begin{proof}  By Theorem 1, $\frac{n-a}{a} \in \mathbb{N}$. Now $$\frac{(n-a)n^3+(n-a-1)n^2+(a-1)n+a}{an^3+(a-1)n^2+(n-a-1)n+(n-a)}$$
$$= \frac{(n-a)(n^3+n^2-n-1)}{a(n^3+n^2-n-1)} = \frac{n-a}{a}.$$
\end{proof}

These $4$-digit solutions were first described by Klosinski and Smolarski~\cite{lK69} in 1969, but their relationship to $5$-digit solutions was not made explicit before now.

It is also interesting to note that $9801$ and $8712$, the two integers in Hardy's discussion of reversals, are included in this family of solutions.

We conclude with the following corollary.

\begin{corollary} There is a $4$-digit solution and a $5$-digit solution for every $n \geq 3$.
\end{corollary}

\begin{proof} Let $a=1$ in the statements of Theorem 1 and Theorem 2 above.
\end{proof}

\subsection*{Some open questions}\label{S:question}

We have shown that there is no $n$ for which there is a $5$-digit solution but no $4$-digit solution.  More specifically, we know that there are $4$- and $5$-digit solutions for every $n \geq 3$.

Although Kaczynski's proof does not generalize directly to $4$- and $5$-digit solutions, it does bring to light several questions about the structure of solutions to the digit reversal problem.

First, it would be interesting to completely characterize $4$- and $5$-digit solutions for $n$.  Namely,

\begin{description}

\item[1.] All known counterexamples to Question 1 occur when $f=0$. Are there counterexamples for which $f \neq 0$? Is there a parameterization for all such counterexamples?

\item[2.] Theorems 1 and 2 exhibit a family of $4$- and $5$-digit solutions for $f=1$ with a particularly nice structure.  To date, no other $4$- or $5$-digit solutions are known for $f=1$. Do such solutions exist?
\end{description}

More generally, 

\begin{description}
\item[3.] Solutions to the digit reversal problem have not been explicitly characterized for more than $5$ digits.  Do there exist analogous results to Theorems 1 and 2 for higher digit solutions?
\end{description}

A Maple package for exploring these questions is available from the author's web page at \texttt{http://www.math.rutgers.edu/\textasciitilde lpudwell/maple.html}.

\subsection*{Acknowledgment}
Thank you to Doron Zeilberger for suggesting this project.

\end{document}